\documentclass[twoside]{article}

\usepackage{amsmath,amsthm,amssymb}

\usepackage{mathptmx}

\usepackage[text={12.5cm,19cm},centering,paperwidth=17cm,paperheight=24cm]{geometry}

\usepackage[fontsize=10.4pt]{scrextend}

\def\titlerunning#1{\gdef\titrun{#1}}

\makeatletter
\def\author#1{\gdef\autrun{\def\and{\unskip, }#1}\gdef\@author{#1}}
\def\address#1{{\def\and{\\\hspace*{15.6pt}}\renewcommand{\thefootnote}{}\footnote{#1}}\markboth{\autrun}{\titrun}}
\makeatother

\def\email#1{email: \href{mailto:#1}{#1} }
\def\subjclass#1{\par\bigskip\noindent\textbf{Mathematics Subject Classification 2020.} #1}
\def\keywords#1{\par\smallskip\noindent\textbf{Keywords.} #1}

\newtheorem{thm}{Theorem}[section]

\theoremstyle{definition}

\numberwithin{equation}{section}

\newcommand{\fconv}{{\mbox{\rm Conv}_{{\rm coe}}(\R^n)}} 
\newcommand{\fconvs}{{\mbox{\rm Conv}_{{\rm sc}}(\R^n)}} 
\newcommand{\fconvfs}{{\mbox{\rm Conv}_{{\rm sc}}(\R^n; \R)}} 
\newcommand{\fconvx}{{\mbox{\rm Conv}(\R^n)}}
\newcommand{\fconvf}{{\mbox{\rm Conv}(\R^n; \R)}}

\newcommand{\sW}{{W^{1,1}}(\R^n)} 
\newcommand{\tW}{W^{1,2}(\R^n)} 

\newcommand{\MA}{\text{\rm MA}} 

\newcommand{\pp}{\mathbin{+_p}}

\newcommand{\cF}{{\mathcal F}}
\newcommand{\cK}{{\mathcal K}} 
\newcommand{\cP}{{\mathcal P}} 
\newcommand{\cQ}{{\mathcal Q}} 
\newcommand\cS{{\mathcal S}}

\newcommand{\sSn}{{\mathcal S}^{n}(\R^n)}

\newcommand{\oM}{\operatorname{M}}
\newcommand{\oid}[1]{\mathop{[ {#1} ]}}
\newcommand{\oD}{\operatorname{D}} 
\newcommand{\fm}{\operatorname{J}} 
\newcommand{\op}{\operatorname{\Pi}} 
\newcommand{\vm}{m} 
\newcommand{\os}[1]{\mathop{\langle {#1} \rangle}} 

\newcommand{\sn}{{\mathbb{S}^{n-1}}}
\newcommand{\Bn}{B^n}

\newcommand{\sln}{\operatorname{SL}(n)}
\newcommand{\gln}{\operatorname{GL}(n)}
\newcommand{\on}{\operatorname{O}(n)}
\newcommand{\son}{\operatorname{SO}(n)}

\newcommand{\Hess}{{\operatorname{D}}^2}

\DeclareMathOperator{\oZ}{\operatorname{Z}}
\DeclareMathOperator{\oY}{\operatorname{Y}}

\newcommand{\oZZ}[2]{\operatorname{V}_{#1,#2}} 

\newcommand{\A}{\mathbb A}

\newcommand{\R}{{\mathbb R}}
\newcommand{\T}{\mathbb T}
\newcommand{\Z}{{\mathbb Z}}

\newcommand{\dom}{\operatorname{dom}}

\newcommand{\dd}{\,\mathrm{d}}

\newcommand{\Had}[2]{D_{#1}^{#2}} 
 
\newcommand{\concave}{\operatorname{\rm Conc}} 

\usepackage[hyperfootnotes=false,colorlinks=true,allcolors=blue]{hyperref}

\begin{document}

\titlerunning{Geometric valuation theory}

\title{\textbf{Geometric valuation theory}}

\author{Monika Ludwig}

\date{}

\maketitle

\address{Institut f\"ur Diskrete Mathematik und Geometrie,
Technische Universit\"at Wien,
Wiedner Hauptstra\ss e 8-10/1046,
1040 Wien, Austria;
\email{monika.ludwig@tuwien.ac.at}}

\begin{abstract}
A brief introduction to geometric valuation theory is given. The focus is on classification results for valuations on convex bodies and on function spaces.
\end{abstract}

\section{Introduction}
Measurement is part of the literal meaning of geometry and geometric valuation theory deals with measurement in the following sense. We want to associate to  a geometric object a real number (or, more generally, an element of an abelian semigroup $\A$). For example, we can associate to a sufficiently regular subset of $\R^n$ its $n$-dimensional volume or the $(n-1)$-dimensional measure of its boundary. Let $\cS$ be a class of subsets of $\R^n$. We call a function $\oZ: \cS\to\A$ a \emph{valuation} if
\begin{equation*}
\oZ(K)+\oZ(L)= \oZ(K\cup L)+\oZ(K\cap L)
\end{equation*}
for all $K,L\in \cS$ with $K\cap L, K\cup L\in\cS$ (and we set $\oZ(\emptyset):=0$). 
Thus, the valuation property is just the inclusion-exclusion principle applied to two sets. In particular, measures on $\R^n$ when restricted to elements of $\cS$ are valuations but there are many additional interesting valuations.

\goodbreak

In his Third Problem, Hilbert asked whether an elementary definition of volume on polytopes is possible. In 1900, it was known that it is possible on $\R^2$ but the question was open in higher dimensions. Let $\cP^n$ be the set of convex polytopes in $\R^n$ and call $\oZ:\cP^n\to \R$ \emph{simple} if $\oZ(P)=0$ for all lower dimensional polytopes. Using our terminology, Hilbert's Third Problem turns out to be equivalent to the question whether every simple, rigid motion invariant valuation $\oZ:\cP^n\to \R$ is a multiple of $n$-dimensional volume for $n\ge3$. Dehn \cite{Dehn1901} solved Hilbert's Third Problem by constructing a simple, rigid motion invariant valuation that is not a multiple of volume and thereby showed that an elementary definition of volume is not possible for $n\ge 3$.

\goodbreak
Blaschke \cite{BlaschkeIntegralH2} took the important next step by asking for classification results for invariant valuations on $\cP^n$ and on the space of convex bodies, $\cK^n$, that is, of non-empty, compact, convex sets in $\R^n$. For a class $\cS$ of subsets of $\R^n$, we say that a function $ \oZ:\cS \to \A$ is \emph{$G$ invariant} for a group $G$ acting on $\R^n$ if
$\oZ(\phi K)= \oZ(K)$ for all $\phi\in G$ and  $K\in \cS$.
Blaschke's question is motivated by Klein's Erlangen Program. 
We will describe some of the results that were obtained in this direction, in particular, focusing on the special linear group, $\sln$, and the group of (orientation preserving) rotations, $\son$. Often additional regularity assumptions are required and for $\A$ a topological semigroup, we consider continuous and upper semicontinuous valuations, where the topology on $\cK^n$ and its subspaces is induced by the Hausdorff metric.

In addition to classification results and their applications, structural results for spaces of valuations have attracted a lot of attention in recent years. We refer to the books and surveys \cite{Alesker_Kent, AleskerFu_Barcelona, Bernig_AIG}. Valuations were also considered on various additional spaces, in particular, on manifolds (see \cite{Alesker07}). We will restrict our attention to subspaces of $\cK^n$ and to recent results on valuations on spaces of real valued functions. On a space $X$ of (extended) real valued functions, a function $\oZ:X\to\A$ is called a \emph{valuation} if
$$\oZ(u)+\oZ(v)=\oZ(u\vee v) + \oZ(u \wedge v)$$
for all $u,v\in X$ such that also their pointwise maximum $u\vee v$ and  pointwise minimum $u\wedge v$ belong to $X$. 
Since spaces of convex bodies can be embedded in various function spaces in such a way that union and intersection of convex bodies correspond to pointwise minimum and maximum of functions, this notion generalizes the classical notion.

\section{Affine valuations on convex bodies}\label{affine}

The first classification result in geometric valuation theory is due to Blaschke. He worked on polytopes and aimed at a complete classification of rigid motion invariant valuations. However, at a certain step, he had to assume also $\sln$ invariance and established the following result (and the corresponding result on polytopes).

\begin{thm}[Blaschke \cite{BlaschkeIntegralH2}]\label{Blaschke}
A functional $\oZ:\cK^n\to \R$ is  a continuous, translation and $\sln$ invariant valuation if and only if  there are  $c_0, c_n\in\R$ such that 
$$\oZ(K) = c_0 V_0(K)+c_nV_n(K)$$
for every $K\in\cK^n$. 
\end{thm}

\noindent
Here, $V_0(K):=1$ is the Euler characteristic of $K$ and $V_n(K)$ its $n$-dimensional volume.
It has become customary to refer to results that involve invariance (or covariance) with respect to $\sln$ as affine results and the title of this section is to be understood in this sense. 

\goodbreak
We will first describe results for affine valuations on polytopes and then on general convex bodies. While on $\cP^n$ a complete classification of  $\sln$ invariant valuations has been established, we require additional regularity assumptions on $\cK^n$. Such assumptions are also used on important subspaces of $\cP^n$  and $\cK^n$. We will also describe results for affine valuations with values in tensor spaces, spaces of convex bodies and related spaces.

\subsection{$\sln$ invariant valuations on convex polytopes}

We call a function $\zeta:[0,\infty)\to \R$ a \emph{Cauchy function} if
$$\zeta(x+y)=\zeta(x)+\zeta(y)$$
for every $x,y\in [0,\infty)$. Cauchy functions are well understood and can be completely described (if we assume the axiom of choice) by their values on a Hamel basis. 

The following result gives a complete classification of translation and $\sln$ invariant valuations on polytopes and is closely related to Theorem \ref{Blaschke}.

\begin{thm}[\cite{LudwigReitzner_AP}]\label{LR_AP}
A functional $\oZ:\cP^n\to \R$ is  a translation and $\sln$ invariant valuation if and only if  there are  $c_0\in\R$  and a Cauchy function $\zeta:[0,\infty)\to \R$ such that 
$$\oZ(P) = c_0 V_0(P)+ \zeta\big(V_n(P)\big)$$
for every $P\in\cP^n$. 
\end{thm}

\noindent
Even without translation invariance, a complete classification can be obtained (see \cite{LudwigReitzner_AP}). Next, we state the case when the valuation is in addition continuous. We write $[0,P]$ for the convex hull of the origin  and $P\in\cP^n$.

\begin{thm}[\cite{LudwigReitzner_AP}]\label{co:3}
A functional $\,\oZ:\cP^n \to \R$ is a continuous and $\,\sln$ invariant valuation if and only if there are  $c_0,  c_n, d_n\in\R$   such that
$$\oZ (P)= c_0 V_0(P) +  c_nV_n(P)+ d_nV_n([0,P])$$
for every $P\in\cP^n$.
\end{thm}

\noindent
Corresponding results are known on the space, $\cP^n_0$, of polytopes containing the origin (see \cite{LudwigReitzner_AP}).

Let $\cP^n_{(0)}$ be the space of convex polytopes  in $\R^n$ that contain the origin in their interiors. Here, we have additional interesting valuations connected to polarity. For $K\in\cK^n$, define its polar by
$$K^*:=\{y\in\R^n: \langle x,y\rangle\le 1\,\, \text{ for all }\,x\in K\},$$
where $\langle x, y\rangle$ is the inner product of $x,y\in\R^n$.
If $P\in\cP^n_{(0)}$, then $P^*\in\cP^n_{(0)}$. Hence, setting
$$V_n^*(P):= V_n(P^*),$$
we obtain a finite valued functional on $\cP^n_{(0)}$ and it follows easily from properties of polarity that it is a valuation.

\goodbreak
Valuations on $\cP_{(0)}^n$ were first considered in  \cite{Ludwig:origin}, where a classification of Borel measurable, $\sln$ invariant and homogeneous valuations was established. Here, we say that $\oZ: \cP_{(0)}^n\to \R$ is \emph{homogeneous} if there is $q\in\R$ such that
$$\oZ(t\,P)=t^q\oZ(P)$$
for every $P\in\cP^n_{(0)}$ and $t>0$. We say that $\oZ$ is \emph{Borel measurable} if the pre-image of every open set is a Borel set. We use corresponding notions on $\cK^n$ and related spaces. 

The results from \cite{Ludwig:origin} were strengthened by Haberl \& Parapatits.

\begin{thm}[Haberl \& Parapatits \cite{Haberl:Parapatits_centro, HaberlParapatits_moments}]\label{HP}
A functional $\oZ:\cP_{(0)}^n\to \R$ is  a Borel measurable and $\sln$ invariant valuation if and only if  there are $c_0, c_n, c_{-n}\in\R$  such that 
$$\oZ(P)= c_0 V_0(P) + c_n V_n(P) + c_{-n} V_n^*(P)$$ 
for every $P\in\cP_{(0)}^n$. 
\end{thm}

\noindent
The regularity assumption is again required to exclude discontinuous solutions of the Cauchy functional equation. It is an open problem to establish a complete classification without such assumption.

We remark that lattice polytopes, that is, convex polytopes with vertices in the integer lattice $\Z^n$, are important in many fields and subjects. The Betke--Kneser theorem \cite{Betke:Kneser} gives a complete classification of valuations on this class that are invariant with respect to translations by integer vectors and by so-called uni\-modular transformations (which can be described by matrices with integer coefficients and determinant $\pm 1$). For more information on valuations on lattice polytopes, see \cite{BoeroeczkyLudwig_survey}.

\subsection{Affine surface areas}

For $K\in\cK^n$, the \emph{affine surface area} of $K$ is defined by 
\begin{equation}\label{asa}
\Omega(K)=\int_{\partial K} \kappa(K,x)^{\frac1{n+1}} \dd x,
\end{equation}
where $\kappa(K,x)$ is the generalized Gaussian curvature of $\partial K$ at $x$ and integration is with respect to the $(n-1)$-dimensional Hausdorff measure.  For smooth convex surfaces, this definition is classical  (see \cite{Blaschke}). It is also classical that $\Omega$ is translation and $\sln$ invariant for smooth surfaces.
The extension of the definition of affine surface area to general convex bodies was obtained more recently in a series of papers by Leichtwei\ss \ \cite{Leichtweiss88}, Lutwak \cite{Lutwak91}, and Sch\"utt \& Werner \cite{Schuett:Werner90}.
There it is also proved that $\Omega$ is translation and $\sln$ invariant on $\cK^n$. 
The notion of affine surface area is fundamental in affine differential geometry. Moreover, since many basic problems in discrete and stochastic geometry are translation and $\sln$ invariant, affine surface area has found numerous applications in these fields  (see  \cite{FejesToth, Gruber}).
It follows easily from \eqref{asa} that $\Omega$ vanishes on polytopes and therefore is  not continuous. The long conjectured upper semicontinuity of affine surface area (for smooth surfaces as well as for general convex surfaces)  was proved  by Lutwak \cite{Lutwak91}. For a proof that $\Omega$ is a valuation, see \cite{Schuett93}. 

\goodbreak
The following result gives a classification of upper semicontinuous, translation and $\sln$ invariant valuations and represents a strengthening of Theorem \ref{Blaschke}. It  provides a characterization of affine surface area.

\begin{thm}[\cite{Ludwig:Reitzner}]
A functional $\oZ:\cK^n\to \R$ is  an upper semicontinuous, trans\-lation and $\sln$ invariant valuation if and only if  there are  $c_0, c_n\in\R$  and $c\ge 0$  such that 
$$\oZ(K)= c_0 V_0(K) +c_n V_n(K)+ c\,\Omega(K)$$
for every $K\in\cK^n$. 
\end{thm}

\noindent
For $n=2$, this result was proved in \cite{Ludwig:affinelength}, where also applications to asymptotic approximation by polytopes were obtained.

A complete classification of translation and $\sln$ invariant valuation on $\cK^n$ appears to be out of reach. Already a weakening of upper semicontinuity to, say, Baire-one (that is, a pointwise limit of continuous functionals) would be interesting and would have applications in discrete and stochastic geometry.

Let $\cK^n_{(0)}$ be the space of convex bodies  in $\R^n$ containing the origin in their interiors. For such a convex body with smooth boundary, the \emph{centro-affine surface area} is a classical notion that can be defined by
\begin{equation*}
\int_{\partial K} \kappa_0(K,x)^{\frac 12} \dd V_K(x),
\end{equation*}
where $\dd V_K(x): =  \langle x, u_K(x)\rangle\dd x$ with $u_K(x)$ the outer unit normal vector to $K$ at $x$ is (up to a constant) the cone measure on $\partial K$ and 
$$\kappa_0(K,x):=\frac{\kappa(K,x)}{\langle x, u_K(x)\rangle^{n+1}}.$$
It is classical that $\Omega_n$ is $\gln$ invariant. 
Lutwak \cite{Lutwak96} extended this notion to general convex bodies in $\cK^n_{(0)}$ and showed that $\Omega_n$ is upper semicontinuous. 

The following result gives a complete classification of upper semicontinuous, $\gln$ invariant valuations on $\cK^n_{(0)}$ and provides a characterization of centro-affine surface area.

\begin{thm}[\cite{Ludwig:Reitzner2}]
A functional $\,\oZ:\cK_{(0)}^n\to \R$ is  an upper semicontinuous, $\gln$ invariant valuation if and only if  there are  $c_0\in\R$  and $c\ge 0$  such that 
$$\oZ(K)= c_0 V_0(K) + c\,\Omega_n(K)$$
for every $K\in\cK_{(0)}^n$. 
\end{thm}

\noindent
Lutwak \cite{Lutwak96} defined so-called $L_p$-affine surface areas which were  characterized in \cite{Ludwig:Reitzner2} as upper semicontinuous, $\sln$ invariant and homogeneous valuations. 

\goodbreak
A more general notion, now called \emph{Orlicz affine surface area}, was introduced in \cite{Ludwig:Reitzner2}. Let
$$\concave[0,\infty):=\Big\{\zeta: [0,\infty)\to [0,\infty): \zeta \text{ concave}, \, \lim_{t\to 0} \zeta(t)=\lim_{t\to\infty} \frac{\zeta(t)}{t}=0\Big\}.$$
The following result gives a classification of upper semicontinuous, $\sln$ invariant valuations on $\cK^n_{(0)}$ and provides a characterization of Orlicz affine surface areas.

\begin{thm}[\cite{Haberl:Parapatits_centro, Ludwig:Reitzner2}]
A functional $\oZ:\cK_{(0)}^n\to \R$ is  an upper semicontinuous, $\sln$~invariant valuation if and only if  there are $c_0, c_n, c_{-n}\in\R$  and $ \zeta\in\concave[0,\infty)$  such that 
$$\oZ(K)= c_0 V_0(K) + c_n V_n(K) + c_{-n} V_n^*(K)+ \int_{\partial K} \zeta(\kappa_0(K,x))\dd V_K(x)$$
for every $K\in\cK_{(0)}^n$. 
\end{thm}

\noindent
Here, the classification of upper semicontinuous, $\sln$ invariant valuations  vanishing on polytopes from \cite{Ludwig:Reitzner2} is combined with Theorem \ref{HP} by Haberl \& Parapatits.

\subsection{Vector and tensor valuations}

We say that $\oZ:\cP^n \to \R^n$ is \emph{$\sln$ equivariant} if
$$
\oZ(\phi P) = \phi \,\oZ(P)$$
for all $\phi \in \sln$ and $P\in\cP^n$. We use corresponding definitions for subspaces of $\cP^n$. 

The study of $\sln$ equivariant  vector valuations on convex polytopes containing the origin in their interiors was started in \cite{Ludwig:moment}, where a  classification of Borel measurable, $\sln$ equivariant, homogeneous valuations was established.  Haberl  \& Parapatits strengthened this result and obtained the following complete classification, which we state for $n\ge 3$. 

\begin{thm}[Haberl \& Parapatits \cite{HaberlParapatits_tensor, HaberlParapatits_moments}]\label{moments}
A function $\,\oZ:\cP^n_{(0)}\to \R^n$ is a Borel measurable and $\,\sln$ equivariant valuation  if and only if there is  $c\in\R$ such that
$$\oZ(P)=c\, \vm(P)$$ 
for every  $P\in\cP^n_{(0)}$. 
\end{thm}

\noindent
Here, for $P\in\cP^n$, the moment vector $\vm(P)$ is defined by
$\vm(P)=\int_P x\dd x$.

\goodbreak
Zeng \& Ma showed that it is possible to obtain a complete classification of vector valuations on convex polytopes without any regularity assumptions. We state their result for $n\ge 3$. 

\begin{thm}[Zeng \& Ma \cite{MaZeng}] A function $\,\oZ:\cP^n\to \R^n$ is an $\sln$ equivariant valuation if and only if there are $c,d\in\R$ such that
$$\oZ(P)=c\, \vm(P)+d\, \vm([0,P])$$
for every $P\in \cP^n$.
\end{thm} 

\noindent
In the same paper, a complete classification result is also established for $n=2$. The obtained valuations depend on Cauchy functions.
\goodbreak

Also higher rank tensor valuations are important in the geometry of convex bodies. In particular, the moment matrix $M^{2,0}(K)$ of a convex body $K$ is a most valuable tool through its connection to the Legendre ellipsoid and the notion of isotropic position. In a certain way dual is the so-called LYZ ellipsoid, which was introduced by Lutwak, Yang \& Zhang \cite{LYZ2000b, LYZ2002}. Associated to this ellipsoid is the LYZ matrix, which was characterized as a matrix valuation on convex polytopes containing the origin  in \cite{Ludwig:matrix}. The LYZ matrix corresponds to the Fisher information matrix  \cite{Ludwig:Fisher, LYZ2000b, LYZ2002} important in statistics and information theory.

Haberl \& Parapatits \cite{HaberlParapatits_tensor} extended the result from \cite{Ludwig:matrix} to general symmetric tensor valuations.  For $p\ge1$, let $\T^p(\R^n)$ denote the space of symmetric $p$-tensor on $\R^n$.  We identify $\R^n$ with its dual space and regard each symmetric $p$-tensor as a symmetric $p$-linear functional on $(\R^n)^p$. We say that $\oZ:\cP^n_{(0)}\to \T^p(\R^n)$ is \emph{$\sln$ equivariant} if 
$$\oZ(\phi P) (y_1,\dots, y_p)=\oZ(P) (\phi^{-1} y_1, \dots, \phi^{-1} y_p)$$
for all $y_1, \dots, y_p\in\R^n$,  all $\phi \in \sln$ and all $P\in\cP^n_{(0)}$. We state the result by Haberl \& Parapatits for $n\ge 3$ and $p\ge 2$.

\goodbreak
\begin{thm}[Haberl \& Parapatits \cite{HaberlParapatits_tensor}]\label{tensor}
A function $\oZ:\cP^n_{(0)}\to \T^p(\R^n)$ is a Borel measurable, $\sln$ equivariant valuation  if and only if there are $c, d\in\R$ such that
$$\oZ(P)=c\, M^{p,0}(P)+ d\, M^{0,p}(P^*)$$
for every $P\in\cP^n_{(0)}$. 
\end{thm}

\noindent Here, the $p$th moment tensor of a convex polytope $P\in\cP^n_{(0)}$ is defined by
\begin{equation}\label{moment_tensor}
M^{p,0}(P):= \frac1{p!}\int_P x^p \dd x
\end{equation}
where $x^p$ is the $p$-fold symmetric tensor product of $x\in\R^n$ and  the $p$th LYZ tensor is
$$M^{0,p}(P) := \int_{\sn} y^p\dd S_{n-1,p}(P, y),$$
where $S_{n-1,p}(P,\cdot)$ is the $L^p$ surface area measure of $P$, which  is a central notion in the $L^p$ Brunn--Minkowski theory (see \cite{Lutwak93b, Lutwak96}). 
\goodbreak

For classifications of matrix valuation on $\cP^n$ without regularity assumptions, see \cite{Ma_matrix, MaWei}, and for tensor valuations on lattice polytopes, see \cite{LudwigSilverstein}. Continuous tensor valuations on complex vector spaces are classified in \cite{ABDK}.

\subsection{Convex body valued valuations and related notions}

Affinely associated convex bodies play an important role in convex geometry (see \cite[Chapter 10]{Schneider:CB2}). We have already mentioned the Legendre and the LYZ ellipsoid and describe here results on valuations $\oZ: \cK^n \to \cK^n$, where we choose suitable additions on $\cK^n$. The most classical choice is the \emph{Minkowski addition}, where for $K,L\in\cK^n$, 
$$K+L:=\{x+y: x\in K, y\in L\},$$
and such valuations are called \emph{Minkowski valuations}.

The first classification result for Minkowski valuations was obtained in \cite{Ludwig:projection} and strengthened in \cite{Ludwig:Minkowski}. It provides a characterization of projection bodies, a notion that was introduced by Minkowski.

\begin{thm}[\cite{Ludwig:Minkowski}]\label{contra_convex}
An operator $\,\oZ:\cP^n\to \cK^n$  is a  translation invariant, $\,\sln$ contra\-variant  Minkowski valuation if and only if  there is $c\ge 0$ such that 
$$\oZ P=c \op P$$
for every  $P\in\cP^n$. 
\end{thm}

\noindent
Here, we describe convex bodies by their support functions, where for $K\in\cK^n$, the \emph{support function} $h(K,\cdot):\R^n\to \R$ is given by
$$h(K,y):=\max\{\langle x,y\rangle: x\in K\}.$$
The support function is homogeneous of degree 1 and convex and any such function is the support function of a convex bodies. For $K\in\cK^n$, the \emph{projection body} of $K$ is defined by
$$h(\op K, y):=V_{n-1}(K\vert y^\perp)$$
for $y\in \sn$ where $y^\perp$ is the hyperplane orthogonal to $y$ and $K\vert y^\perp$ denotes the image of the orthogonal projection of $K$ onto $y^\perp$. We say that $\oZ: \cP^n\to \cK^n$ is \emph{$\sln$ contravariant} if
$$\oZ(\phi P)= \phi^{-t} \oZ P$$
for all $\phi\in\sln$ and $P\in\cP^n$, where $\phi^{-t}$ is the inverse of the transpose of $\phi$. For more information on projection bodies and their many applications, see \cite{Gardner, Schneider:CB2}.

We say that $\oZ: \cP^n\to \cK^n$ is \emph{$\sln$ equivariant} if
$$\oZ(\phi P)= \phi \oZ P$$
for all $\phi\in\sln$ and $P\in\cP^n$. 
The following result establishes a classification $\sln$ equi\-variant valuations. 

\begin{thm}[\cite{Ludwig:Minkowski}]\label{equi_convex}
An operator $\,\oZ:\cP^n\to \cK^n$  is a  translation invariant, $\,\sln$ equi\-variant  Minkowski valuation if and only if  there is $c\ge 0$ such that 
$$\oZ P=c \oD P$$
for every $P\in\cP^n$. 
\end{thm}

\noindent
Here, the operator $P \mapsto \oD P =\{x-y: x,y\in P\}$ assigns to $P$ its {\em difference body} (see \cite{Gardner, Schneider:CB2}).
\goodbreak

A classification of $\sln$ equivariant, homogeneous Minkowski valuations on the space, $\cK^n_0$, of convex bodies containing the origin was obtained  in \cite{Ludwig:Minkowski}. The result was strengthened by Haberl \cite{Haberl_sln}, who was able to drop the assumption of homogeneity.
Let $n\ge 3$.

\begin{thm}[Haberl \cite{Haberl_sln}]\label{JEMS2}
An operator $\,\oZ :\cK^n_0\to \cK^n$ is a continuous, 
$\,\sln\!$ equivariant Minkowski valuation if and only if 
there are $c_0\in\R$ and $c_1,c_2, c_3\ge 0$ such that
$$\oZ K = 
c_0\, \vm(K) + c_1\, K +c_2(- K) +c_3 \oM K $$
for every $K\in\cK^n_0$. 
\end{thm}

\noindent
Here, the  {\em moment body}, $\,\oM K$, of $K$ is defined by 
\begin{equation*}
h(\oM K, y):= \int_ K \vert\langle x, y\rangle\vert\dd x
\end{equation*}
for $y\in\R^n$.
When divided by the volume of $K$, the moment body of $K$ is called its \emph{centroid body} and is a classical and important notion going back to at least Dupin (see \cite{Gardner, Schneider:CB2}). Results corresponding to Theorem \ref{JEMS2} for $\sln$ contravariant Minkowski valuations were obtained in  \cite{Haberl_sln, Ludwig:Minkowski}. 
On the space, $\cP_0^n$, of convex polytopes containing the origin classification results for $\sln$ contravariant Minkowski valuations were established  in \cite{Haberl_sln, Ludwig:Minkowski} without assuming continuity and additional operators appear. For the $\sln$ equivariant case, such results were established in \cite{LiLeng_polytopes}.

We remark that the results from Theorem \ref{JEMS2} and the corresponding results in the $\sln$ equivariant case were complemented in \cite{Schuster:Wannerer, Wannerer2011} by classification results for continuous, homogeneous Minkowski valuations on $\cK^n$. A complete classification for $\sln$ equivariant Minkowski valuations on $\cP_0^n$ was established in \cite{Haberl_sln}.
On the space of convex bodies that contain the origin in their interiors, moment bodies allow to define $\sln$ equivariant Minkowski valuations using polarity. For continuous, $\sln$ equivariant, homogeneous valuations,  a complete classification on this space was established in \cite{Ludwig:convex}. For Minkowski valuations on lattice polytopes, see \cite{BoeroeczkyLudwig}.

Classification results for Minkowski valuations on complex vector spaces were established by Abardia \& Bernig \cite{Abardia, Abardia15, Abardia:Bernig}. They introduce and characterize complex projection and difference bodies.

\goodbreak
An important extension of the classical Brunn--Minkowski theory is the more recent 
$L^p$ Brunn--Minkowski theory (see \cite{Lutwak93b, Lutwak96}). For $p>1$, the $L_p$ sum of $K,L\in\cK_0^n$ is defined by
$$h^p(K \pp L,y):= h^p(K,y)+ h^p(L,y)$$
for $y\in\R^n$.  An $L^p$ Minkowski valuations $\oZ: \cK^n\to \cK_0^n$ is a valuation where on $\cK_0^n$ this addition is chosen. Classification results were obtained in \cite{LiLeng_polytopes, Ludwig:Minkowski, Parapatits:co, Parapatits:contravariant} and led to the definition of asymmetric $L^p$ projection and moment bodies (see \cite{Ludwig:Minkowski}). Inequalities for these new classes of operators were established by Haberl \& Schuster \cite{Haberl:Schuster1}. They generalize the $L^p$ Petty projection  and the $L^p$ Busemann--Petty moment inequalities, which were established by Lutwak, Yang \& Zhang \cite{LYZ2000}, and were, in turn, generalized within the Orlicz--Brunn--Minkowski inequality by Lutwak, Yang  \& Zhang \cite{LYZ2010a, LYZ2010b}. For information on valuations in this setting, see \cite{LiLeng_orlicz}.

A classical notion of addition on full dimensional convex bodies in $\R^n$ is  Blaschke addition, which is defined using the sum of surface area measures of convex bodies and the solution of the classical Minkowski problem. So-called \emph{Blaschke valuations} were classified in \cite{Haberl:blaschke}. For the corresponding question within the $L^p$ Brunn--Minkowski theory, see \cite{LiYuanLeng}.

The dual Brunn--Minkowski theory, established by Lutwak \cite{Lutwak75}, is, in a certain way, dual to the classical theory. Star bodies replace convex bodies and radial addition (defined by the addition of radial functions) corresponds to Minkowski addition. Intersection bodies in the dual Brunn--Minkowski theory correspond to projection bodies in the classical theory. Intersection bodies were critical in the solution of the Busemann--Petty problem \cite{Lutwak88, Zhang99a}. A classification of radial valuations and a characterization of the intersection body operator was established in
\cite{Ludwig:intersection}. Replacing radial addition by $L^p$~radial addition leads to $L^p$ radial valuations (see \cite{Haberl:star, Haberl:Ludwig} for classification results).

Since convex bodies can be described by support functions and star bodies by radial functions, a natural extension of the results described above is a classification of valuations $\oZ: \cK^n\to F(\R^n)$, where $F(\R^n)$ is a suitable space of functions on $\R^n$. Such results were obtained by Li \cite{Li2020, Li2021} and by Li \& Ma \cite{LiMa} where a characterization of the Laplace transform on convex bodies is established. Another way to describe convex bodies is by suitable measures and a  classification of measure valued valuations was obtained by Haberl \& Parapatits \cite{Haberl:Parapatits_crelle}, where characterization results of surface area measures and of $L^p$ surface area measures were established.

\section{The Hadwiger theorem on convex bodies}

The classical Steiner formula states that the volume of the outer parallel set of a convex body at distance $r>0$  can be expressed as a polynomial in $r$ of degree at most $n$.  Using that the outer parallel set of  $K\in\cK^n$ at distance $r>0$ is just the Minkowski sum of $K$ and $r\Bn$, the ball of radius $r$, we get
\begin{equation*}
V_n({K+r \Bn})=\sum_{j=0}^n r^{n-j}\kappa_{n-j} V_j(K)
\end{equation*}
for every $r>0$, where $\kappa_j$ is the $j$-dimensional volume of the unit ball in $\R^j$ (with the convention that $\kappa_0:=1$). The coefficients $V_j(K)$ are known as the \emph{intrinsic volumes} of $K$. Up to  normalization and numbering, they coincide with the classical quermass\-integrals. In particular,  $V_{n-1}(K)$ is  proportional to the surface area of $K$ and $V_1(K)$ to its mean width (cf. \cite{Schneider:CB2}). 

\goodbreak
The celebrated Hadwiger theorem gives a characterization of intrinsic volumes and a complete classification of continuous, translation and rotation invariant valuations.
For $n=2$, it follows from the positive solution to Hilbert's Third Problem in this case. It was proved for $n=3$ in \cite{Hadwiger51} and then for general $n\ge3$ in \cite{Hadwiger52}.

\begin{thm}[Hadwiger \cite{Hadwiger52}]\label{hugo}
A functional $\oZ:\cK^n\to \R$ is  a continuous, translation and rotation invariant valuation if and only if  there are  $c_0, \ldots, c_n\in\R$ such that 
$$\oZ(K) = c_0  V_0(K)+\dots+c_n V_n(K)$$
for every $K\in\cK^n$. 
\end{thm}

\noindent
The Hadwiger theorem leads to effortless proofs of numerous results in integral geo\-metry and geometric probability (see \cite{Hadwiger:V,Klain:Rota}).
An alternate proof of the Hadwiger theorem is due to Klain \cite{Klain95}.

We will describe results on translation invariant and rotation equivariant valuations with values in tensor spaces and spaces of convex bodies. We remark that upper semicontinuous, translation and rotation invariant valuations were only classified in the planar case (see \cite{Ludwig:semi}).

\subsection{Vector and tensor valuation}

The first classification of vector valuations was
 established by Hadwiger  \& Schneider \cite{Hadwiger:Schneider} using rotation equivariant valuations $\oZ:\cK^n\to \R^n$, that is,  valuations such that
$$\oZ(\phi K)= \phi \oZ(K)$$
for all  $\phi\in\son$ and $K\in\cK^n$.

\begin{thm}[Hadwiger \& Schneider \cite{Hadwiger:Schneider}]\label{HadwigerSchneider}
A function $\,\oZ:\cK^n\to \R^n$ is  a continuous, translation covariant, rotation equivariant valuation if and only if  there are $c_1,\ldots, c_{n+1}\in\R$ such that 
$$
\oZ(K) =  c_1 \oM^{1,0}_1(K)+\dots+c_{n+1}\oM^{1,0}_{n+1}(K)$$
for every $K\in\cK^n$. 
\end{thm}

\noindent
Here $\oM^{1,0}_i(K):=\Phi^{1,0}_i(K)$ are the \emph{intrinsic vectors} of $K$ (see (\ref{steiner_f2}) below) and  see (\ref{covariant}) for the definition of translation covariance.

\goodbreak
The theorem by Hadwiger \& Schneider was extended by Alesker \cite{Alesker98, Alesker00a} (based on \cite{Alesker99}) to  a classification of continuous, translation covariant, rotation equi\-variant tensor valuations on $\cK^n$. Just as the intrinsic volumes can be obtained from the Steiner polynomial, the moment tensor (defined in \eqref{moment_tensor})
satisfies the Steiner formula 
\begin{equation}\label{steiner_f2}
\oM^{p, 0}(K+r \,B^n)= \sum_{j=0}^{n+p} r^{n+p-j} \kappa_{n+p-j} \,\sum_{k\ge 0} \Phi^{p-k,k}_{j-p+k}(K)
\end{equation}
for $K\in\cK^n$ and $r\ge0$. 
The coefficients $\Phi_k^{p,s}(K)$ are called the \emph{Minkowski tensors} of $K$ (see \cite[Section~5.4]{Schneider:CB2}).
Recall that $\T^p(\R^n)$ is the space of symmetric $p$-tensors on $\R^n$ and let $Q\in\T^2(\R^n)$ be the metric tensor, that is, $Q(x,y)= \langle x, y\rangle$ for $x,y\in\R^n$. 

\begin{thm}[Alesker \cite{Alesker98}]\label{Alesker}
A function $\,\oZ:\cK^n\to \T^p(\R^n)$ is  a continuous, trans\-lation covariant, rotation equivariant valuation if and only if  $\,\oZ$  can be written as linear combination of the symmetric tensor products
$Q^l \,\Phi_k^{m,s}$ with  $2l+m+s=p$.
\end{thm}

\goodbreak
\noindent
Here, a valuation $\oZ: \cK^n \to
\T^p(\R^n)$ is called \emph{translation covariant} if there exist associated
functions $\oZ^j: \cK^n \to \T^j(\R^n)$ for $j=0,\dots, p$
such that
\begin{equation}\label{covariant}
  \oZ(K+y) =\sum_{j=0}^p \oZ^{r-j}(K) \frac{y^{j}}{j!}
\end{equation}
for all $y\in\R^n$ and $K\in\cK^n$, where on the right side we sum over symmetric tensor products.  We say that $\oZ$ is  \emph{$G$ equivariant} for a group $G$ acting on $\R^n$ if
 \begin{equation*}
  \oZ(\phi K) (y_1, \dots, y_p) = \oZ(K) (\phi^t y_1, \dots, \phi^t y_p)
\end{equation*}
for all $y_1, \dots,y_p\in\R^n$,  all transformation $\phi\in G$, and all $K\in\cK^n$, where $\phi^t$ is the transpose of $\phi$.

For a classification of local tensor valuations, see \cite{HugSchneider14}, and for applications in various fields, including astronomy and material sciences, see \cite{KiderlenVedelJensen}.

\subsection{Convex body valued valuations}

An operator $\oZ: \cK^n\to\cK^n$ is called \emph{Minkowski additive} if
$$\oZ(K+L)= \oZ(K)+ \oZ(L)$$
for all $K, L\in \cK^n$. Since $K+L= K\cup L+ K\cap L$ for $K, L\in\cK^n$ with $K\cup L\in\cK^n$, it is easy to see that every Minkowski additive operator is a Minkowski valuation. While the first classification results for Minkowski valuations were established in \cite{Ludwig:projection}, Schneider  \cite{Schneider74} earlier obtained the first classification results for rotation equi\-variant Minkowski additive operators. For continuous, translation invariant, rotation equi\-variant Minkowski valuations, so far no complete classification has been established. But the following representation is known to hold. Let $\mathcal{M}_{\mathrm{cen}}\left(\mathbb{S}^{n-1}\right)$ and
$C_{\mathrm{cen}}\left(\mathrm{S}^{n-1}\right)$ denote the spaces of signed Borel measures and continuous functions on $\mathrm{S}^{n-1},$ respectively, having their center of mass at the origin.

\begin{thm}[Schuster \& Wannerer \cite{SchusterWannererJems}]\label{Schuster:Wannerer}
If $\,\oZ: \mathcal{K}^{n} \rightarrow \mathcal{K}^{n}$ is a continuous, translation invariant, rotation equivariant Minkowski valuation, then there are uniquely determined constants $c_{0}, c_{n} \geq 0$ and $\,\mathrm{SO}(n-1)$ invariant measures $\mu_{i} \in \mathcal{M}_{\mathrm{cen}}\left(\mathbb{S}^{n-1}\right)$ for $1 \leq i \leq n-2,$  as well as an $\mathrm{SO}(n-1)$
invariant function $\zeta_{n-1} \in C_{\mathrm{cen}}\left(\mathrm{S}^{n-1}\right)$ such that
$$
h({\oZ K}, \cdot)=c_{0}+\sum_{i=1}^{n-2} S_{i}(K, \cdot) * \mu_{i}+S_{n-1}(K, \cdot) * \zeta_{n-1}+c_{n} V_{n}(K)
$$
for every $K \in \mathcal{K}^{n}$.
\end{thm}

\noindent
The Borel measures $S_{i}(K, \cdot)$ on $\mathbb{S}^{n-1}$ are Aleksandrov's area measures (see \cite{Schneider:CB2}) of $K \in \mathcal{K}^{n}$. The convolution of functions and measures on $\mathbb{S}^{n-1}$ is induced from the group $\son$ by identifying $\mathbb{S}^{n-1}$ with the homogeneous space $\mathrm{SO}(n) / \mathrm{SO}(n-1)$ (see \cite{SchusterWannererJems}). The above representation formula has to be read in the sense of equality of measures and $h({\oZ K}, \cdot)$ is identified with the measure with this density.

\section{More on invariant valuations on convex bodies}

Translation invariant valuations on polytopes were classified using simplicity or mild regularity assumptions. 
Hadwiger \cite{Hadwiger52a} established a complete classification of simple, weakly continuous, translation invariant valuations on convex polytopes. Here, informally, a valuation is \emph{weakly continuous} if it is continuous under parallel displacements
of the facets of a polytope. Hadwiger's result was extended by McMullen \cite{McMullen1983} to the following result. 

\begin{thm}[McMullen \cite{McMullen1983}]
A functional $\oZ: \cP^n\to \R$ is a weakly continuous,
translation invariant valuation if and only if
$$\oZ(P)=\sum_{j=0}^n \sum_{\,F\in \cF_j(P)} \oY_j(N(P,F)) \,V_j(F)
$$
for every $P\in\cP^n$
where $\oY_j: \cQ^{n-j}\to\R$ is a simple valuation.
\end{thm}

\noindent
Here, $\cF_j(P)$ is the set of $j$-dimensional faces of $P$ and $N(P,F)$ is the normal cone to $P$ at $F$ while $\cQ^k$  is the system of all closed polyhedral convex cones of dimension at most $k$.  We remark that valuations on convex polyhedral cones (or, equivalently, on spherical polytopes) are not yet well understood and the problems to classify simple, rotation invariant valuations on spherical polytopes and on spherical convex bodies are open on spheres of dimension $\ge 3$ (even if continuity is assumed).
Kusejko \& Parapatits \cite{KusejkoParapatits} extended Hadwiger's result  and established a complete classification of simple, translation invariant valuations on polytopes using Cauchy functions. 

\goodbreak
Hadwiger \cite{Hadwiger:V}  proved that simple, continuous, translation invariant valuations on $\cK^n$ have a \emph{homogeneous decomposition}. His result was extended by McMullen \cite{McMullen77}.

\begin{thm}[McMullen \cite{McMullen77}]\label{McM} If $\,\oZ: \cK^n\to \R$ is a continuous and translation invariant valuation, then 
$$\oZ=\oZ_0+\dots +\oZ_n$$
where $\oZ_j: \cK^n \to \R$ is a continuous, translation invariant valuation
that is homo\-geneous of degree  $j$.
\end{thm}

\noindent
It is easy to see that every continuous, translation invariant valuation that is homo\-geneous of degree $0$ is a multiple of the Euler characteristic. For the degrees of homogeneity $j=n$ and $j=n-1$, the following results hold.

\begin{thm}[Hadwiger \cite{Hadwiger:V}]\label{hugo_n}
A functional $\oZ:\cP^n\to \R$ is  a translation invariant valuation that is homogeneous of degree $n$ if and only if  there is $c\in\R$ such that 
$$\oZ(P) =c\, V_n(P)$$
for every $P\in\cP^n$. 
\end{thm}

\goodbreak
\begin{thm}[McMullen \cite{McMullen80}]\label{mcmullen_n-1}
A functional $\oZ:\cK^n\to \R$ is  a continuous, trans\-lation invariant valuation which is homogeneous of degree $(n-1)$ if and only if  there is $\zeta\in C(\sn)$ such that 
$$\oZ(K) =\int_{\sn} \zeta(y)\dd S_{n-1}(K,y)$$
for every $K\in\cK^n$. The function $\zeta$ is uniquely determined up to addition of the restriction of a linear function.
\end{thm}

\noindent
Continuous, translation invariant valuations that are homogeneous of degree $1$ were classified by Goodey \& Weil  \cite{GoodeyWeil1984}. 

While a complete classification of continuous, translation invariant valuations on $\cK^n$ is out of reach, Alesker \cite{Alesker01} proved the following result.

\begin{thm}[Alesker \cite{Alesker01}]\label{dense}
For $\,0\le j \le n$, linear combinations of the valuations
$$\big\{K\mapsto V(K[j], K_1, \dots, K_{n-j}):  K_1, \dots, K_{n-j}\in\cK^n\big\}$$
are dense in the space of continuous, translation invariant valuations that are homo\-geneous of degree $j$. 
\end{thm}

\noindent
Here, $V(K[j], K_1, \dots, K_{n-j})$ is the mixed volume of $K$ taken $j$ times and $K_1, \dots, K_{n-j}$ while the topology on the space of continuous, translation invariant valuations is induced by the norm $\Vert \oZ \Vert:= \sup\{\vert \oZ(K)\vert: K\in\cK^n,\, K\subseteq B^n\}$. Alesker's result confirms a conjecture by McMullen \cite{McMullen80} and is based on Alesker's so-called irreducibility theorem, which was proved in \cite{Alesker01} and which has far-reaching consequences.

\goodbreak

For simple valuations,  the following complete classification was established by Klain \& Schneider.

\begin{thm}[Klain \cite{Klain95}, Schneider \cite{Schneider:simple}]\label{klain_schneider}
A functional $\oZ:\cK^n\to \R$ is  a simple, continuous, translation invariant valuation if and only if  there are $c\in \R$ and an odd function $\zeta\in C(\sn)$ such that 
$$\oZ(K) =\int_{\sn} \zeta(y)\dd S_{n-1}(K,y)+c\, V_n(K)$$
for every $K\in\cK^n$. The function $\zeta$ is uniquely determined up to addition of the restriction of a linear function.
\end{thm}

\noindent
Klain \cite{Klain95} used his classification of simple valuations in his proof of the Hadwiger theorem. For an alternate proof of Theorem \ref{klain_schneider}, see \cite{KusejkoParapatits}.

A valuation $\oZ: \cK^n\to \R$ is called \emph{translatively polynomial} if $x\mapsto \oZ(P+x)$ is a polynomial in the coordinates of $x\in\R^n$ for all $K\in\cK^n$. 
Alesker  \cite{Alesker99} established a complete classification of continuous, translatively polynomial, rotation invariant valuations on $\cK^n$. Theorem \ref{Alesker} is the version of this result for tensor valuations.

\goodbreak
Classification results for continuous, translation invariant valuations that are invariant under indefinite orthogonal groups were established by Alesker \& Faifman  \cite{AleskerFaifman} and Bernig \& Faifman \cite{BernigFaifman2017}. For subgroups of the orthogonal group $\on$, the following result holds.

\begin{thm}[Alesker \cite{Alesker00, Alesker07}]
For a compact subgroup $G$ of \,$\on$, the linear space of continuous, translation and $G$ invariant valuations on $\cK^n$ is finite dimensional if and only if
$G$ acts transitively on $\sn$.
\end{thm}

\noindent
As the classification of the such subgroups $G$ is known, it was a natural task (which was already proposed in \cite{Alesker00}) to find bases for spaces of $G$ invariant valuations  (see \cite{Alesker01, Alesker2004, Alesker05, Alesker2008, Bernig2009, Bernig2011, Bernig2012, Bernig:Fu, BernigSolanes2014, BernigSolanes2017, BernigVoide} for results on real valued valuations and \cite{BDS21, Wannerer14} for results on tensor and measure valued valuations).

\section{Affine valuations on function spaces}

We describe classification results for valuations on function spaces that correspond to the results in Section \ref{affine}. Let $F(\R^n)$ be a space of functions $f: \R^n\to [-\infty, \infty]$ and let $G$ be a subgroup of $\gln$. An operator $\oZ: F(\R^n)\to \A$ is \emph{$G$ invariant} if  
$$\oZ(f\circ\phi^{-1})=\oZ(f)$$
for all $\phi\in G$ and $f\in F(\R^n)$.  If $G$ acts on $\A$, we say that an operator $\oZ: F(\R^n)\to \A$ is \emph{$G$ contravariant}
if for some $q\in\R$,
$$\oZ(f\circ \phi^{-1})=\vert\det \phi \vert^q \,\phi^{-t} \oZ (f)$$
for all $\phi\in G$ and $f\in F(\R^n)$. It is \emph{$G$ equivariant} if for some $q\in\R$,
$$\oZ(f\circ \phi^{-1})=\vert\det \phi \vert^q \,\phi \oZ (f)$$
for all $\phi\in G$ and $f\in F(\R^n)$.
It is called \emph{homogeneous} if for some $q\in\R$,
$$\oZ(s\, f)= \vert s\vert^q \oZ(f)$$
for all $s\in\R$ and $f\in F(\R^n)$ such that $s\,f\in F(\R^n)$. An operator is called \emph{affinely contravariant} if it is trans\-lation invariant, $\gln$ contravariant and homogeneous.

\subsection{Valuations on Sobolev spaces}\label{Sobolev}

For $p\ge 1$, let $W^{1,p}(\R^n)$ be the Sobolev space of functions belonging to $L^p(\R^n)$ whose distributional first-order derivatives belong to $L^p(\R^n)$.

The following result corresponds to Theorem \ref{contra_convex}. Let $\cK^n_c$ be the set of origin-symmetric convex bodies in $\R^n$. Let $n\ge 3$.

\begin{thm}[\cite{Ludwig:SobVal}]\label{mink}
An operator $\,\oZ :\sW\to \cK_c^n$ is a  continuous, affinely contra\-variant Minkowski
valuation 
if and only if  
there is  $c \ge0$   such that
$$\oZ (f) = c \op \os{f}$$
for every $f\in\sW$. 
\end{thm}

\noindent
Here, for $f\in\sW$, the \emph{LYZ body} $\os{f}$  is defined by Lutwak, Yang \& Zhang \cite{LYZ2006} as the unique origin-symmetric convex body in $\R^n$ such that 
\begin{equation}\label{LYZbody}
\int_{\sn} \zeta(y)\dd S_{n-1}(\os{f},y)=\int_{\R^n} \zeta(\nabla f (x))\dd x
\end{equation}
for every even continuous function $\zeta:\R^n\to \R$ that is homo\-geneous of degree~1. 
Equation~\eqref{LYZbody}  is a functional version of the classical even Minkowski problem. 

Combined with \eqref{LYZbody}, it follows from the definition of projection bodies and surface area measures that for $f\in\sW$ and $y\in\sn$,
$$
h(\op \os{ f}, y)=\frac12\int_{\R^n} \vert\langle \nabla f(x),y\rangle\vert\dd x.
$$
We remark that the convex body $\os{f}\,$ has proved to be critical in geometric analysis: the affine Sobolev--Zhang inequality \cite{Zhang99} is a volume inequality for the polar body of $\op \os{f}$, which strengthens and implies the Euclidean case of the classical  Sobolev inequality, and it was proved in \cite{LYZ2006} that $\os{f}$ describes the optimal Sobolev norm of $f\in \sW$.
Tuo Wang\ \cite{Tuo_Wang} studied the LYZ operator $f\mapsto \os{f}$ on the space of functions of bounded variation. Here, the LYZ operator is not a valuation anymore but Wang~\cite{Tuo_Wang_semi} established a characterization as an affinely covariant Blaschke semi-valuation.

\goodbreak
The following classification of tensor valuations corresponds to Theorem \ref{tensor} for $p=2$. Let $n\ge 3$.

\begin{thm}[\cite{Ludwig:Fisher}]\label{fisher}
An operator $\,\oZ :\tW\to \T^2(\R^n)$ is a continuous, affinely contra\-variant valuation if and only if  
there is   $c \in\R$   such that
$$\oZ (f) = c \fm(f^2)
$$
for every $f\in\tW$. 
\end{thm}

\noindent
Here,  we write $\fm(h)$ for the \emph{Fisher information matrix} of the weakly differentiable function $h:\R^n\to[0,\infty)$, that  is, the $n\times n$ matrix with entries
\begin{equation}\label{fisherdef}
\fm_{ij}(g) =\int_{\R^n} \frac{\partial \log h(x)}{\partial x_i}\,\frac{\partial \log h(x)}{\partial x_j}\, h(x)\dd x.
\end{equation}
We remark that the Fisher information matrix plays an important role in information theory and statistics (see \cite{Cover:Thomas}). In general, Fisher information is a measure of the minimum error in the maximum likelihood estimate of  a parameter in a distribution. The Fisher information matrix \eqref{fisherdef} describes such an error  for a random vector of density $h$ with respect to a location parameter. 

For results on real valued valuations on Sobolev spaces, see \cite{Ma2016}.

\subsection{Valuations on convex functions} 

We write $\fconvx$ for  the space of convex functions $u:\R^n \to (-\infty, \infty]$ that are lower semicontinuous and proper, that is, $u \not\equiv\infty$. We equip $\fconvx$ and its subspaces with the topology induced by epi-convergence (see \cite{RockafellarWets}). Let
$$\fconv:=\{u\in\fconvx: \lim\nolimits_{\vert x\vert \to \infty} u(x)=\infty\}$$
be the space of \emph{coercive}, convex functions, where $\vert x\vert$ is the Euclidean norm of $x\in\R^n$. The following result corresponds to Theorem \ref{Blaschke}.

\begin{thm}[\cite{Colesanti-Ludwig-Mussnig-1}]
A functional $\,\oZ :\fconv\to [0,\infty)$ is a continuous, trans\-lation and $\sln\!$ invariant valuation if and only if there are a continuous function $\zeta_0: \R\to[0,\infty)$  and a continuous function $\zeta_n: \R\to[0,\infty)$  with finite $(n-1)$th moment  such that 
\begin{equation*}
\oZ (u) = \zeta_0\big(\min\nolimits_{x\in\R^n}u(x)\big) + \int_{\dom u} \zeta_n\big(u(x)\big)\dd x
\end{equation*}
for every $u\in\fconv$. 
\end{thm}
\noindent
Here, a function $\zeta:\R \to [0,\infty)$ has finite $k$th moment if $\int_0^{+\infty} t^{k} \zeta(t)\dd t<+\infty$ and $\dom u$ is the \emph{domain} of $u$, that is, $\dom u=\{x\in\R^n: u(x)<+\infty\}$. 

\goodbreak
Let $\fconvf$ be the space of finite valued convex functions, that is, of convex functions $u: \R^n\to \R$. We say that $u\in \fconvx$ is \emph{super-coercive} if  
$ \lim\nolimits_{\vert x\vert \to \infty} \frac{u(x)}{\vert x\vert}=\infty$. Let 
$\fconvfs$ be the space of \emph{super-coercive}, finite valued, convex functions.
The following result corresponds to Theorem \ref{HP}.

\begin{thm}[Mussnig \cite{Mussnig21}]
A functional $\,\oZ :\fconvfs\to [0,\infty)$ is a con\-tinuous, translation and $\sln\!$ invariant valuation if and only if there are a continuous $\zeta_0:\R\to[0,\infty)$, a continuous  $\zeta_n:\R\to[0,\infty)$  with finite $(n-1)$th moment, and a continuous  $\zeta_{-n}:\R\to[0,\infty)$  whose support is bounded from above such that 
\begin{equation*}
\oZ (u) = \zeta_0\big(\min\nolimits_{x\in\R^n}u(x)\big) + \int_{\R^n} \zeta_n\big(u(x)\big)\dd x + \int_{\R^n} \zeta_{-n} (u(x)) \dd \MA(u, x)
\end{equation*}
for every $u\in\fconvfs$. 
\end{thm}
\noindent
Here, $\MA(u, \cdot)$ denotes the Monge--Amp\`ere measure of $u$, which is also called the $n$th Hessian measure.
See  \cite{Mussnig19} for a result on coercive functions in $\fconvf$.

\goodbreak

The following results correspond to Theorem \ref{contra_convex} and Theorem \ref{equi_convex}. Let $n\ge 3$.

\begin{thm}[\cite{Colesanti-Ludwig-Mussnig-2}]\label{contravariant}
An operator $\,\oZ :\fconv\to \cK^n$ is a  continuous, monotone,  translation invariant, $\sln\!$ contra\-variant Minkowski valuation if and only if there is a continuous, decreasing  $\zeta: \R\to[0,\infty)$ with finite $(n-2)$th moment such that
 $$\oZ (u) = \op \os{\zeta \circ u}$$
for every $u\in\fconv$. 
\end{thm}

\noindent
For $u\in\fconv$ and suitable $\zeta\in C(\R)$, define the \emph{level set body} $\oid{\zeta \circ u}\in\cK^n$ by
$$
h(\oid{\zeta \circ u}, y)=  \int_0^{+\infty} h(\{\zeta\circ u\geq t\},y) \dd t
$$
for $y\in\R^n$. Hence the level set body is a Minkowski average of the level sets.

\begin{thm}[\cite{Colesanti-Ludwig-Mussnig-2}]\label{Mcovariant}
An operator $\,\oZ :\fconv\to \cK^n$ is a continuous, monotone,  translation invariant, $\sln$ equivariant  Minkowski valuation if and only if there is a continuous, decreasing  $\zeta:\R\to[0,\infty)$ with finite integral over $[0,\infty)$ such that
 $$\oZ (u) = \oD \oid{\zeta \circ u}$$
for every $u\in\fconv$. 
\end{thm}

\noindent
We remark that the results in this section can easily be translated to classification results for valuations on log-concave functions. In this setting, the results on convex body valued valuations were strengthened by Mussnig \cite{Mussnig_log}.

\section{The Hadwiger theorem on convex functions}

We call a functional $\oZ: \fconvs\to \R$ \emph{epi-translation invariant} if
$$\oZ(u\circ \tau^{-1}+c)= \oZ(u)$$
for all translations $\tau: \R^n\to \R^n$ and $c\in\R$. Hence $\oZ(u)$ is not changed by translations of the epi-graph of $u$. To state the Hadwiger theorem on $\fconvf$, we need to define functional versions of the intrinsic volumes.
Let $C_b((0,\infty))$ be the set of continuous functions on $(0,\infty)$ with bounded support. For $0\leq j \leq n-1$, let
$$
\Had{j}{n}:=\big\{\zeta\in C_b((0,\infty))\colon  \lim_{s\to 0^+} s^{n-j} \zeta(s)=0,  \lim_{s\to 0^+} \int_s^{\infty}  t^{n-j-1}\zeta(t) \dd t \text{ exists and is finite}\big\}.
$$
In addition, let $\Had{n}{n}$ be the set of functions $\zeta\in C_b((0,\infty))$ where $\lim_{s\to 0^+} \zeta(s)$ exists and is finite and set $\zeta(0):=\lim_{s\to 0^+} \zeta(s)$.

\begin{thm}[\cite{Colesanti-Ludwig-Mussnig-5}]\label{main one way} 
	For $0\le j\le n$ and $\zeta\in\Had{j}{n}$, 
	there exists a unique, continuous, epi-translation and 	rotation invariant valuation 
	$\oZZ{j}{\zeta}\colon\fconvs\to\R$ such that 
	\begin{equation*}
	\oZZ{j}{\zeta}(u)=  \int_{\R^n} \zeta(\vert\nabla u(x)\vert) [\Hess u(x)]_{n-j} \dd x
	\end{equation*}
	for every $u\in\fconvs\cap C^2_+(\R^n)$. 
\end{thm}

\noindent
Here, $\Hess u$ is the Hessian matrix of $u$ and $[\Hess u(x)]_k$  the $k$th elementary symmetric functions of the eigenvalues of $\Hess u(x)$ (with the convention that $[\Hess u(x)]_0:\equiv 1$) while $C^2_+(\R^n)$ is the space of twice continuously differentiable functions with positive definite Hessians.
We remark that $\oZZ{0}{\zeta}$ is  constant on $\fconvs$.

The following result is the Hadwiger theorem on $\fconvs$. Let $n\ge 2$.

\begin{thm}[\cite{Colesanti-Ludwig-Mussnig-5}]
\label{thm:hadwiger_convex_functions}
A functional $\,\oZ:\fconvs \to \R$ is a continuous, epi-translation and rotation invariant valuation if and only if there are functions $\zeta_0\in\Had{0}{n}$, \dots, $\zeta_n\in\Had{n}{n}$  such that
\begin{equation*}
\oZ(u)= \oZZ{0}{\zeta_0}(u) +\dots +\oZZ{n}{\zeta_n}(u)
\end{equation*}
for every $u\in\fconvs$.
\end{thm}

\noindent 
A comparison of Theorem \ref{hugo} and Theorem \ref{thm:hadwiger_convex_functions} shows that for $0\le j\le n$ and $\zeta\in \Had{j}{n}$, the functional $\oZZ{j}{\zeta}$ plays a role corresponding to that of the $j$th intrinsic volume $V_j$. Hence, we call $\oZZ{j}{\zeta}$ a $j$th \emph{functional intrinsic volume} on $\fconvs$.

\goodbreak
We call a functional $\oZ: \fconvf\to \R$ \emph{dually epi-translation invariant} if
$$\oZ(v + \ell +c)= \oZ(v)$$
for all linear functions $\ell: \R^n\to \R$ and $c\in\R$. Using the convex conjugate or Legendre transform of $u\in\fconvs$, given by
$$u^*(y):=\sup\nolimits_{x\in\R^n} \big(\langle x,y \rangle - u(x) \big)$$
for $y\in\R^n$, we see that $v\mapsto \oZ(v)$ is dually epi-translation invariant on $\fconvf$ if and only if $u\mapsto \oZ(u^*)$ is epi-translation invariant on $\fconvs$.
It was proved in \cite{Colesanti-Ludwig-Mussnig-3} that $\oZ$ is a continuous valuation on $\fconvf$ if and only if $\oZ^*\colon\fconvs\to\R$,  defined by
$$
\oZ^*(u):=\oZ(u^*),
$$
 is a continuous valuation on $\fconvs$. This fact permits us to transfer results  valid for valuations on $\fconvs$ to results for valuations on 
$\fconvf$ and vice versa.  

\goodbreak
The following result is obtained from Theorem \ref{main one way} by using convex conjugation.

\begin{thm}[\cite{Colesanti-Ludwig-Mussnig-5}]\label{dual main one way} 
For $0\le j\le n$ and $\zeta\in\Had{j}{n}$, 
the functional
$\oZZ{j}{\zeta}^*\colon\fconvf\to\R$ is a continuous, dually epi-translation and 	rotation invariant valuation  such that
\begin{equation}\label{dhessian}
	\oZZ{j}{\zeta}^*(v)=  \int_{\R^n} \zeta(\vert x\vert) [\Hess v(x)]_{j} \dd x
\end{equation}
for every $v\in\fconvf\cap C^2_+(\R^n)$. 
\end{thm}

\noindent
Here, $\oZZ{j}{\zeta}^*(v):= \oZZ{j}{\zeta}(v^*)$ for $0\le j\le n$ and $\zeta\in\Had{j}{n}$.
Theorem \ref{thm:hadwiger_convex_functions} has the following dual version. Let $n\ge 2$.

\begin{thm}[\cite{Colesanti-Ludwig-Mussnig-5}]
\label{dthm:hadwiger_convex_functions}
A functional $\oZ:\fconvf \to \R$ is a continuous, dually epi-translation and rotation in\-va\-riant valuation if and only if there are  functions $\zeta_0\in\Had{0}{n}$, \dots, $\zeta_n\in\Had{n}{n}$  such that
\begin{equation*}
\oZ(v)= \oZZ{0}{\zeta_0}^*(v) +\dots +\oZZ{n}{\zeta_n}^*(v)
\end{equation*}
for every $v\in\fconvf$.
\end{thm}

Applications of the Hadwiger theorem on convex functions including integral geometric formulas and additional representations of functional intrinsic volumes will be presented in follow-up papers to \cite{Colesanti-Ludwig-Mussnig-5}.

\section{More on invariant valuations on function spaces}

For continuous, epi-translation invariant valuations on $\fconvs$,  the existence of a homogeneous decomposition corresponding to Theorem \ref{McM} was established in \cite{Colesanti-Ludwig-Mussnig-4}, that is, every such valuation is a linear combination of continuous, epi-translation invariant valuations that are epi-homogeneous of degree $j$ and $0\le j\le n$. Here $\oZ$ is called \emph{epi-homogeneous} of degree $j$ if $\oZ(u)$ is multiplied by $t^j$ when the epi-graph of $u$ is multiplied by $t>0$. It is not difficult to see that every continuous, epi-translation invariant valuation that is epi-homogeneous of degree $0$ is constant. 

\goodbreak
The following classification corresponding to Theorem \ref{hugo_n} was established in \cite{Colesanti-Ludwig-Mussnig-4}.

\begin{thm}[\cite{Colesanti-Ludwig-Mussnig-4}]\label{clm_n}
A functional $\oZ:\fconvs\to \R$ is  an epi-translation invariant valuation that is epi-homogeneous of degree $n$ if and only if  there is $\zeta\in C_c(\R^n)$ such that 
$$\oZ(u) =\int_{\dom u} \zeta(\nabla u(x)) \dd x$$
for every $u\in\fconvs$. 
\end{thm}

\noindent
Here, $C_c(\R^n)$ is the space of continuous functions with compact support.
The result corresponding to Theorem \ref{clm_n} on $\fconvf$ is stated next.

\begin{thm}[\cite{Colesanti-Ludwig-Mussnig-4}]
A functional $\oZ:\fconvf\to \R$ is  a dually epi-translation invariant valuation that is homogeneous of degree $n$ if and only if  there is $\zeta\in C_c(\R^n)$ such that 
$$\oZ(v) =\int_{\R^n} \zeta(x) \dd \MA(v,x)$$
for every $v\in\fconvf$. 
\end{thm}

\noindent
See \cite{Colesanti-Ludwig-Mussnig-4}, for more information on homogeneous decompositions and why such results do not hold for many spaces of convex functions.
 For more results on valuations on convex functions, see  \cite{Alesker_cf, CavallinaColesanti, Knoerr1, Knoerr2}, and for results on valuations on quasi-concave functions, see \cite{ColesantiLombardi, ColesantiLombardiParapatits}.

While formally not results for valuations on function spaces, classification results for valuations on star shaped sets in $\R^n$ were the motivation for some of the results on function spaces. Let $\sSn$ be the space of sets $S\subset \R^n$ which are star shaped with respect to the origin and whose radial functions $\rho(S, \cdot): \sn\to [0,\infty]$, given by
$$\rho(S,x):= \sup\{r\ge 0: r\,x \in S\},$$ 
are in $L^n(\sn)$. Let 
$\cS_0$ be the space of star bodies, that is, of star shaped sets with continuous radial functions. We remark that $\cS^n_0$ is the space used in the dual Brunn--Minkowski theory (see \cite{Gardner,Lutwak75}). 
Note that union and intersection on $\sSn$ and on $\cS^n_0$ correspond to the pointwise maximum and minimum for radial functions. We equip $\sSn$ with the topology induced by the $L^n$ norm on $\sn$  and $\cS^n_0$ with the topology induced by the maximum norm.

\goodbreak
Klain \cite{Klain97} established the following classification results on star shaped sets.

\begin{thm}[Klain \cite{Klain97}]\label{klain}
A functional $\,\oZ : \sSn\to\R$ is a continuous, rotation invariant valuation  with  $\oZ(\{0\})=0$ 
if and only if there is $\zeta\in C([0,\infty))$ with the properties that $\zeta(0)=0$ and
$\vert\zeta(t)\vert \le c+ d \vert t\vert^n$ for all $\,t\in\R$\, for some $c, d\ge0$ such that
$$\oZ(S)=\int_{\sn} \zeta(\rho(S, y ))\dd y$$
for every $S \in \sSn$.
\end{thm}

\noindent If the valuation $\oZ$ in Theorem \ref{klain} is in addition positively homogeneous of degree $p$, then $\zeta(t) =c\, t^p$ with $c\in \R$ and $0\le p\le n$ and hence $\oZ$ is a dual mixed volume (as defined by Lutwak  \cite{Lutwak75}). 

\goodbreak
Tsang \cite{Tsang:Lp} obtained classification results for valuations on $L^p(X,\mu)$, when $X$ is a non-atomic measure space. Here we state a special case of his results that complements Theorem \ref{klain}. Let $p\ge 1$.

\begin{thm}[Tsang \cite{Tsang:Lp}]
A functional $\oZ : L^p(\R^n)\to\R$ is a continuous, translation invariant valuation that vanishes on the null function
if and only if there is  $\zeta\in C(\R)$ with the property 
 that $\vert\zeta(t)\vert \le c \,\vert t\vert^p$ for all $\,t\in\R$ for some $c\ge 0$ such that
$$\oZ(f)=\int_{\R^n} \zeta(f(x)) \dd x$$
for every $f \in L^p(\R^n)$.
\end{thm}

\noindent
We remark that also Theorem \ref{klain} can be written as a classification result on the space of non-negative functions in $L^n(\sn)$ (also see \cite{Tsang:Lp}). For results on tensor and Minkowski valuations on $L^p$ space, see \cite{Ludwig:MM, Ober2014,Tsang:Minkowski}.

\goodbreak
Villanueva \cite{Villanueva2016} obtained classification results for non-negative valuations on star bodies. In \cite{TradaceteVillanueva2017}, Tradacete \& Villanueva showed that a result corresponding to the classification from Theorem \ref{klain} is valid on $\cS^n_0$. A complete classification on $\cS^n_0$ is given in the following result.

\goodbreak

\begin{thm}[Tradacete \&Villanueva \cite{TradaceteVillanueva}]\label{TV}
A functional $\,\oZ : \cS^n_0\to \R$ is a continuous valuation 
if and only if there are a finite Borel measure $\mu$ on $\sn$ and a function $\zeta:[0,\infty)\times \sn\to \R$ that fulfills the strong Carath\'eodory condition with respect to $\mu$  such that
$$\oZ(S)=\int_{\sn} \zeta(\rho(S, y ), y)\dd \mu(y)$$
for every $u \in \cS^n_0$.
\end{thm}

\noindent
Here, we say that $\zeta:[0,\infty)\times \sn\to \R$ fulfills the \emph{strong  Carath\'eodory condition} with respect to $\mu$ if $\zeta(s, \cdot)$ is Borel measurable for all $s\ge 0$ and $\zeta(\cdot, y)$ is continuous for $\mu$ almost every $y\in\sn$ while for every $t>0$ there is $\xi_t\in L^1(\sn, \mu)$ such that $\zeta(s, y)\le \xi_t(y)$ for $s<t$ and $\mu$ almost every $y\in\sn$. We remark that Theorem \ref{TV} can be rewritten as a result on valuations on non-negative functions in $C(\sn)$.

Classification results for valuations on Lipschitz functions on $\sn$ were obtained in
\cite{ColesantiPagniniTradaceteVillanueva, ColesantiPagniniTradaceteVillanueva2}  
 and on Banach lattices in \cite{TradaceteVillanueva_imrn}. A Hadwiger theorem for valuations on definable functions was established in \cite{BaryshnikovGhristWright}.

\goodbreak
\small

\end{document}